\documentclass{article}
\usepackage{graphicx} 
\usepackage{setspace}
\usepackage[a4paper,left=2.5cm,right=2.5cm,top=2.5cm,bottom=2.5cm]{geometry}
\usepackage{amsmath, amssymb, amsthm}
\usepackage{mdframed}

\title{A Survey of Alternative Solutions to the Congruum Problem}
\author{Nathanael Johnson}
\date{13 April 2025}

\begin{document}

\begin{doublespace}

\maketitle

\begin{abstract}
    A congruum was first defined by Leonardo Pisano in 1225 and it is defined as the common difference in an arithmetic progression of three perfect squares. Later that year in his book \emph{Liber Quadratorum}, Pisano proved that congruums can never perfect squares themselves, a finding that was later revisited by Pierre de Fermat in 1670. His proof is now known as Fermat's Right Triangle Theorem. In this paper, four alternative proofs to Pisano's original proof are demonstrated and offered with each proof requiring a different scope of mathematical knowledge. The proofs are by direct Diophantine analysis, parameterization of differences, Heronian triangle construction, and infinite descent. In showing these proofs, it is demonstrated that there are alternatives to the method of decomposing perfect squares as sums of odd numbers as Pisano did in his proof in 1225.
\end{abstract}

\section{Introduction}

In 1225, Leonardo Pisano, commonly known as Fibonacci, published his \emph{Liber Quadratorum} (Book of Squares), in which he laid the groundwork for many concepts of number theory, algebra, and arithmetic progressions. In particular, Proposition 14, sometimes referred to as the Congruum Problem or Congruum Theorem, asserts and proves the statement that an arithmetic progression of three perfect squares cannot have a common difference (called a congruum when the common difference is between perfect squares) that is also a perfect square. Equivalently, there are no four integers a, b, c, d such that $a^2+b^2=c^2$ and $b^2+c^2=d^2$. In this representation $b^2$ represents the congruum. However, Fibonacci's proof is several pages long and requires the representation of each perfect square, $a^2$, $c^2$, and $d^2$, as their corresponding sums of odd numbers. These sums are then manipulated and evaluated into different cases and it is found that the value of any congruum will be of the form $4mn(m+n)(m-n)$ where m and n are coprime integers \cite{fibonacci}. This form is then proven to never equal a perfect square. The simpler result that $mn(m+n)(m-n)$ is never a perfect square is proven again by Pierre de Fermat in 1670, and has become known as Fermat's Right Triangle Theorem \cite{fermat}. In this paper, three alternative proofs for Fibonacci's Congruum Theorem are shown, each independent of each other, each independent of Fibonacci's original proof, and each completed without the need for substituting the perfect squares as corresponding sums of odd numbers as did Fibonacci in his original proof.

\section{Preliminary Information}

\noindent First, a few previously proven and well-known lemmas should be recalled and established as they will be helpful for the purposes of the proof. Since these lemmas are all previously proven and in no way original work of this paper, their proofs will be cited in the work cited and summarized briefly but will not be fully worked out in the text of the main body of the paper.

\newmdtheoremenv{lemma}{Lemma}

\begin{lemma}[Fermat's Right Triangle Theorem \cite{fermat}]
    The area of a right triangle with three integer side lengths is never a perfect square.
\end{lemma}

\noindent The proof is by infinite descent. Essentially, Fermat starts by showing the area of any right triangle is equal to $mn(m^2-n^2)$ (which can quickly be seen using Euclid's Formulae as well) and then by using the properties of m and n being coprime, he is able to construct a smaller triangle with an area equal to a perfect square. Since this cannot continue infinitely for natural numbers, the original area cannot be a perfect square. \cite{fermat}

\begin{lemma}[Euclid's Formulae \cite{euclid}]
    Suppose a, b, c are a Pythagorean triple. Then, using two coprime integers of different parity m,n, and a scaling factor k where k=1 for primitive Pythagorean triples, the sides a, b, and c can parameterized as follows:
    
~~~~~~~~~~~~~~~ $a=k(m^2-n^2)$ ~~~ $b=2kmn$ ~~~ $c=k(m^2+n^2)$ 
\end{lemma}

\noindent The parameterization begins by taking the equation and rearranging it to obtain $\frac{b}{c-a}=\frac{c+a}{b}$ and then setting both fractions equal to $\frac{m}{n}$ for m,n coprime. From here, he was able to solve for $\frac{c}{b}$ and $\frac{a}{b}$ and use properties of m and n to obtain the formulae for primitive a, b, c. Non-primitive triples are then generated by multiplying the corresponding primitive triple by the scaling factor k. \cite{euclid}

\begin{lemma}[Heron's Triangle Formula \cite{heron}]
    For any triangle with side lengths of a, b, c and a semiperimeter s where $2s=a+b+c$, the formula for the area of the triangle A is given by:

~~~~~~~~~~~~~~~    $A=(s(s-a)(s-b)(s-c))^{\frac{1}{2}}$ 
\end{lemma}

\noindent There are many ways to prove this formula but one way would be to split the triangle into two right triangles with the height as the shared side and the base as the side split between the two right triangles. It is then possible to find h in terms of the sides of each right triangle and combine these formulae to obtain the height in terms of a, b, c. The area of the triangle is then one half times the base (which is either a, b, or c) times the height and proper simplification and manipulation will result in Heron's formula.

\begin{lemma}
    For any primitive Pythagorean triple a, b, c, c must be odd.
\end{lemma}

\noindent This is not attributed as any single person's work so it is worth working through the entire proof. Lemma 4 is simply because if a is odd and b is even or vice versa, c must be odd since the square of an odd is an odd and the square of an even is an even, and an odd plus an even is an odd. If a and b are both even, c must be even for the same reason but then all three are even so the triple is no longer primitive. Finally, if a and b are both odd, all odd squares are of the form $(2k+1)^2=4k^2+4k+1$ so they have a remainder of 1 mod 4. Adding two such numbers yields that $c^2$ should be congruent to 2 mod 4. But an even number squared is congruent to 0 mod 4 so this is impossible. The only case left is that c is odd in all primitive cases.

\begin{lemma}
    Suppose a, b, c, d are integers so that $a^2+b^2+c^2=d^2$. Then, using four integers m, n, p, q such that $m+n+p+q \in 2\mathbb{N}+1$, the four integers a, b, c, and d can be parameterized as follows:
    
~~~ $a=m^2+n^2-p^2-q^2$ ~~~ $b=2(mq+np)$ ~~~ $c=2(nq-mp)$ ~~~ $d=m^2+n^2+p^2+q^2$ 
\end{lemma}

This proof requires more nuanced analysis but the general methods are similar to those for Pythagorean triples. The main difference is that several substitutions have to be made throughout to gradually convert the equation into an equation that is easier to analyze. If any of these lemmas requires deeper study to be certain of comprehension, this one is so.

\section{Proof by Direct Diophantine Analysis}

\noindent With these simple lemmas established, it is time to go through the alternative proofs to Fibonacci's Congruum Theorem. The first proof is by direct Diophantine analysis of the resultant system of equations.

\begin{proof}
\noindent We begin with the fact that $a^2+b^2=c^2$ and $b^2+c^2=d^2$. Suppose there did exist a, b, c, d satisfying this system of Diophantine equations such that a, b, c, d are all integers. Then, we substitute the first equation into the second to obtain $a^2 + 2b^2=d^2$ which can be rearranged to equal $2b^2=(d-a)(d+a)$. Then,

$\frac{2b}{d-a} = \frac{d+a}{b} = \frac{m}{n} ~\exists m,n \in \mathbb{N}, ~gcd(m,n)=1$

\noindent This parameterization then  allows each fraction to be compared to $\frac{m}{n}$ so that we have the system of equations:

$\frac{d}{b}+\frac{a}{b}=\frac{m}{n}$

$\frac{d}{b}-\frac{a}{b}=\frac{2n}{m}$

\noindent Therefore, solving for $\frac{d}{b}$ and $\frac{a}{b}$, we obtain:

$\frac{d}{b}=\frac{m^2+2n^2}{2mn}$

$\frac{a}{b}=\frac{m^2-2n^2}{2mn}$

\noindent Alternatively, by beginning with $\frac{b}{d-a} = \frac{d+a}{2b} = \frac{m}{n} ~\exists m,n \in \mathbb{N}, ~gcd(m,n)=1$, the comparable results are:

$\frac{d}{b}=\frac{2m^2+n^2}{2mn}$

$\frac{a}{b}=\frac{2m^2-n^2}{2mn}$

\noindent Since m and n are coprime, they are not both even. If m is odd, the numerator is odd so 2 does not divide the numerator. Since m and n are coprime, neither m nor n divides the numerator either. So the fractions in the first solution are in simplest form. Similarly, if n is odd, the fractions in the second solution are in simplest form. So every such triple a, b, d is represented by at least one of these two results as follows:

Solution 1: $a=2m^2-n^2, ~b=2mn, ~d=2m^2+n^2$

Solution 2: $a=m^2-2n^2, ~b=2mn, ~d=m^2+2n^2$

\noindent Now, we evaluate $a^2+b^2=c^2$ in each case. We obtain the following simplified equations:

Solution 1: $4m^4+n^4=c^2$

Solution 2: $m^4+4n^4=c^2$

\noindent In the case of solution 1, we have a new Pythagorean triple produced by $2m^2, ~n^2, ~c$ which yields an area of $m^2n^2$ for the triangle defined by this new triple. Similarly, in the case of solution 2, we have a new Pythagorean triple produced by $m^2, ~2n^2, ~c$ which yields an area of $m^2n^2$ for the triangle defined by this new triple.

\noindent In both cases, the resultant area is a perfect square which by Lemma 1 is not possible! This a contradiction. Thus, no such solution of four integers a, b, c, d exists for the given equations.
\end{proof}

\noindent As a side note, it is worth mentioning that the derivation of the solution forms of the triples a, b, d is the same general process that can be utilized for deriving Lemma 2 by simply replacing the numerator or denominator of 2b with another copy of b instead.

\section{Proof by Parameterization of Differences}

\noindent Continuing, the second proof is given by directly parameterizing the resultant arithmetic progression. This proof is most similar to the original proof by Fibonacci in that it finds a similar looking form for all congruums. However, the method for arriving at this form is nonetheless very different.

\begin{proof}

\noindent Suppose not. Suppose instead there is a solution such that a, b, c, d are all integers. We begin by noting that $a^2, ~c^2, ~d^2$ form an arithmetic progression with common difference r. We wish to show that r cannot be a perfect square so that $r=b^2$ does not allow for b to be an integer.

\noindent Since $d^2>c^2>a^2$, we know that $d>c>a$ so $c=a+km$ and $d=c+kn=a+km+kn$ where we take m and n to be coprime and $k=gcd(d-c,c-a)$. Then,

$r=d^2-c^2=(a+km+kn)^2-(a+km)^2=k^2n^2+2k^2mn+2akn$

$r=c^2-a^2=(a+km)^2-a^2=k^2m^2+2akm$

\noindent Dividing each formula for r through by k, setting them equal, and solving for a, we have:

$kn^2+2kmn+2an=km^2+2am$

$a=\frac{kn^2-km^2+2kmn}{2(m-n)}=\frac{kmn}{m-n}-\frac{k}{2}(m+n)$

\noindent Now since $a^2+b^2=c^2$ and $b^2+c^2=d^2$, if a, b, c do not form a primitive triple, the common factor will divide $b^2+c^2$ which is $d^2$ so the common factor also divides the other triple. For this reason, we can assume the two triples are both primitive, since if they are not, the scaling factor can be divided out of all four of a, b, c, d. Thus, since both triples are primitive, c and d are odd by Lemma 4 since they are primitive hypotenuses. So b is even since $d^2-c^2=b^2$. So a is odd since $c^2-b^2=a^2$. Thus, d-c and c-a are both even and so $k=gcd(d-c, c-a)$ is even. So the second term of $\frac{k}{2}(m+n)$ is an integer. So $\frac{kmn}{m-n}$ is an integer. But m and n are coprime so both are coprime to m-n. Thus, $J=\frac{k}{m-n}$ is an integer and we can substitute this value in.

\noindent Returning to our value for a, we can simplify now by putting J into the single fraction formula for a and then adding km to obtain c and kn to obtain d. We have:

$a=\frac{J}{2}(n^2-m^2+2mn)$

$c=\frac{J}{2}(m^2+n^2)$

$d=\frac{J}{2}(m^2-n^2+2mn)$

\noindent Then, $r=d^2-c^2$=$\frac{J^2}{4}((m^2-n^2+2mn)^2-(m^2+n^2)^2)=\frac{J^2}{4}(4mn(m^2-n^2))=J^2mn(m^2-n^2)$. So for r to be a perfect square, $J^2mn(m^2-n^2)$ has to be a perfect square meaning $mn(m^2-n^2)$ has to be a perfect square. But this represents the area of the right triangle with side lengths $2mn, ~m^2-n^2, ~m^2+n^2$ (not necessarily primitive, this triple can be confirmed by plugging the values into the Pythagorean equation) which means this right triangle must have a perfect square area which is impossible by Lemma 1! This is a contradiction. Thus, r is not a perfect square and $r=b^2$ can result in b being an integer for such integers a, c, d.
\end{proof}

\noindent As previously stated, this proof is unique among the three in that in remotely resembles the results of Fibonacci's proof method. However, even the result of $J^2mn(m^2-n^2)$ differs is this purely because of the difference in the actual process of the proof. However, this proof circumvents the need to represents squares as sums of odds unlike Fibonacci's proof.

\section{Proof by Heronian Triangle Construction}

\noindent Moving forward, the third proof is given by parameterizing the triples a, b, c and b, c, d using Lemma 3, and then constructing a new right triangle from the parameters. Just as in Proof 2, this proof assumes both triples are primitive and the steps are deliberately omitted since they are previously explained above (though the location of the steps is cited).

\begin{proof}

\noindent Suppose not. Suppose there did exist a, b, c, d that satisfy the system of equations and are all integers. As previously shown in Proof 2, it can be assumed that a, b, c and b, c, d are both primitive triples and b must be even. Then, by Lemma 2 and putting k=1, we have the following parameterizations for each triple:

$a=m^2-n^2$ ~ $b=2mn$ ~ $c=m^2+n^2$

$c=u^2-v^2$ ~~~ $b=2uv$ ~~~ $d=u^2+v^2$

\noindent By transitivity of equality, we have that $mn=uv$ and $m^2+n^2=u^2-v^2$ or $u^2=m^2+n^2+v^2$. 

\noindent Consider the triangle with side lengths $x=u^2-m^2, ~y=u^2-n^2, ~z=u^2-v^2$. By Lemma 3, we have the area for the triangle as follows for $2s=a+b+c$ so that $s=\frac{u^2-m^2+u^2-n^2+u^2-v^2}{2}=\frac{3u^2-u^2}{2}=\frac{2u^2}{2}=u^2$. We obtain:

$A=(u^2(u^2-u^2+m^2)(u^2-u^2+n^2)(u^2-u^2+v^2))^\frac{1}{2}$

$A=(u^2m^2n^2v^2)^\frac{1}{2}=uvmn=uv(uv)=u^2v^2$

\noindent So the area of this triangle, since $mn=uv$, is a perfect square. Thus this triangle cannot be a right triangle. It is easy to check whether $x^2+y^2=z^2$ as follows:

$x^2+y^2=(u^2-m^2)^2+(u^2-n^2)^2$

$x^2+y^2=u^4+m^4-2u^2m^2+u^4+n^4-2u^2n^2=m^4+n^4+2u^2(u^2-m^2-n^2)$

\noindent But since $u^2=m^2+n^2+v^2$ and $mn=uv$, we have the next few steps using these two substitutions respectively resulting in the following:

$x^2+y^2=m^4+n^4+2u^2v^2=m^4+n^4+2m^2n^2=(m^2+n^2)^2=(u^2-v^2)^2=z^2$

\noindent So the triangle is indeed a right triangle with a perfect square area which is impossible! This is a contradiction. So there is no parameterization for both a, b, c and b, c, d being Pythagorean triples.
\end{proof}

\noindent It is important to note that Proof 3 can also be completed by calculating the area of the constructed triangle by the normal base times height method. However, this requires the triangle to be a right triangle which had not yet been established.

\section{Proof by Infinite Descent}

\noindent Finally, the fourth proof is by directly demonstrating infinite descent from any stipulated solution. Like Proofs 2 and 3, this proof uses the fact that a, b, c and b, c, d are both primitive triples. Also, the parameterizations for a, b, c, d are repeated as in Proof 3.

\begin{proof}

\noindent Suppose not. Suppose there existed two triples a, b, c and b, c, d that result from integers a, b, c, d. Then, we parameterize a, b, c and b, c, d as primitive Pythagorean triples using Lemma 2 as below:

$a=m^2-n^2$ ~ $b=2mn$ ~ $c=m^2+n^2$

$c=u^2-v^2$ ~~~ $b=2uv$ ~~~ $d=u^2+v^2$

\noindent As with Proof 3, by transitivity of equality, we have that $mn=uv$ and $m^2+n^2=u^2-v^2$ or $u^2=m^2+n^2+v^2$. In particular, $mn=uv$ can be further examined by representing each side of the equation as two groups of prime factors. each copy of each prime factor, not necessarily distinct, is either in m or n on the LHS and either in u or v on the RHS. Let w, x, y, z be the products of the copies that are in m on the LHS and u on the RHS, m on the LHS and v on the RHS, n on the LHS and u on the RHS, and n on the LHS and v on the RHS, respectively. We then can represent this as follows:

$mn=rs=wxyz$

$m=wx, ~n=yz, ~r=wy, ~s=xz$

\noindent Notice, since m, n are coprime and r, s are coprime, all four of w, x, y, z are pairwise coprime. Also, b should be greater than 0 so all four are also greater than 0. Since $m^2+n^2=u^2-v^2$, we substitute all four variables in terms of w, x, y, z and obtain: 

$w^2x^2+y^2z^2=w^2y^2-x^2z^2$

\noindent Now since $r^2=m^2+n^2+s^2$, and r and s are coprime, we have a primitive Pythagorean quadruple. So then by Lemma 5, we have:

$r=i^2+j^2+k^2+l^2$, $i+j+k+l \in 2\mathbb{N}+1$

\noindent So either one or three of the four i, j, k, l are odd. So the same is true for their squares. So r is also odd. Since r and s are parameters for a primitive triple, they are different parities so s is even. Also since r is odd and $r=wy$, it is known that w and y are each odd. Also one of m and n are odd so either x is odd or z is odd but not both since $s=xz$ is even. Since flipping x and z does not change s and only changes m and n up to symmetry of switching w and y, we can choose x to be the even one of the two and recognize that the exact same procedure from here forward can be completely comparably with x odd and z even.

\noindent If x is even, $x^2$ and $y^2$ are different parity. We can rearrange the equation to obtain the following:

$w^2(y^2-x^2)=z^2(y^2+x^2)$

\noindent To be clear from before, if we had presumed x to odd, we would rearrange so to have $x^2(w^2+z^2)=y^2(w^2-z^2)$ and follow the upcoming steps comparably. Continuing with $w^2(y^2-x^2)=z^2(y^2+x^2)$, we can divide through in two different ways to obtain the equations:

$w^2=\frac{z^2(y^2+x^2)}{y^2-x^2}$

$y^2+x^2=\frac{w^2(y^2-x^2)}{z^2}$

\noindent We see both fractions are equation to integers so they are also integers. Suppose a prime divides both $x^2+y^2$ and $x^2-y^2$. Then, that prime divides $2x^2$ an $2y^2$ since would divide the sum and difference of both. Since x and y are coprime, the primes that divides both must be 2. But x and y are different parity so both $x^2+y^2$ and $x^2-y^2$ are odd. Thus, there is no such prime and they are coprime as well. So $y^2-x^2$ divides $z^2$ since $w^2=\frac{z^2(y^2+x^2)}{y^2-x^2}$ is an integer and $y^2-x^2$ does not divide $x^2+y^2$. Similarly, w, x, y, z are pairwise coprime so since $y^2+x^2=\frac{w^2(y^2-x^2)}{z^2}$ is an integer and $z^2$ does not divide $w^2$, we know $z^2$ divides $y^2-x^2$. So $z^2=y^2-x^2$ since they divide each other and canceling (since $z^2>0$), $w^2=y^2+x^2$. In summary, we have:

$x^2+z^2=y^2$

$x^2 +y^2 =w^2$

\noindent So we have produced another set of four as was stipulated about a, b, c, d where x is the analog of b. But w, x, y, z > 0 so $0<x\leq wxyz<2wxyz=b$. So, if x is even and z is odd, we have produced from a given set of a, b, c, d another set of four w, x, y, z where the common difference $x^2$ is less than the common difference $b^2$ but still greater than 0. Similarly, if x is odd and z is even, we have produced from a given set of a, b, c, d another set of four w, x, y, z where the common difference $z^2$ is less than the common difference $b^2$ but still greater than 0. So in all cases this result is obtained. This would then be possible forever which cannot be so for natural numbers! Thus, by infinite descent, there is no such set of four integers a, b, c, d.
\end{proof}

This proof is unique in that it does not require the use of Fermat's Right Triangle Theorem. It only requires the elementary properties of coprime numbers and it produces a direct link between the theorem and the infinite descent being brought into play.

\section{Conclusion}

\noindent Overall, this has a survey of alternative methods of proof for Fibonacci's Proposition 14 of \emph{Liber Quadratorum}. More specifically, the statement is that there is no arithmetic progression of three perfect squares with a common difference that is also a perfect square. The resultant examination crossed fields of algebra, number theory, and geometry and it is important to note the scope that one proposition can have in the broader proof scene. Future work would most likely relate to using the four methods, Direct Diophantine Analysis, Parameterization of Differences, Heronian Triangle Construction, and Infinite Descent for more complicated problems in the field of Diophantine equations and systems of Diophantine equations. As has already been noted, the method of Proof 1 is completely analogous to the usual way of finding Euclid's formulae for the Pythagorean equation and Proof 2 gives a more direct way of performing the same line thought as did Fibonacci in his original proof.

\pagebreak

\end{doublespace}

\end{document}